\documentclass[11pt]{amsart}
\usepackage{amscd}
\usepackage{amsfonts,amssymb,latexsym}  
\input xy
\xyoption{all}
\xyoption{arc}
\xyoption{ps}
\xyoption{dvips}
\CompileMatrices

\newcommand{\SH}{{\mathcal{H}}}
\newcommand{\SI}{{\mathcal{I}}}

\newcommand{\SM}{{\mathcal{M}}}

\newcommand{\SO}{{\mathcal{O}}}

\renewcommand{\SS}{{\mathcal{S}}}

\newcommand{\SU}{{\mathcal{U}}}

\renewcommand{\AA}{\mathbb{A}}
\newcommand{\CC}{\mathbb{C}}

\newcommand{\tr}{{\rm tr}}
\newcommand{\Lie}{\operatorname{Lie}}

\newcommand{\Aut}{\operatorname{Aut}}

\newcommand{\isom}{\cong}
\newcommand{\Ext}{\operatorname{Ext}}

\newcommand{\Hom}{\operatorname{Hom}}

\newcommand{\Pic}{\operatorname{Pic}}
\newcommand{\Spec}{\operatorname{Spec}}
\newcommand{\Sym}{\operatorname{Sym}}
\newcommand{\id}{\operatorname{id}}

\newcommand{\surj}{\twoheadrightarrow}

\newcommand{\too}{\longrightarrow}
\newcommand{\rk}{\operatorname{rk}}

\newcommand{\wt}{\widetilde}

\newtheorem{proposition}{Proposition}[section]
\newtheorem{theorem}[proposition]{Theorem}

\newtheorem{lemma}[proposition]{Lemma}

\numberwithin{equation}{section}

\title[Torelli theorem for Higgs bundles]{A Torelli
theorem for the moduli space of Higgs bundles
on a curve}
\author[I. Biswas and T. L. G\'omez]{Indranil Biswas and
Tom\'as L. G\'omez}
\date{13 July 2001}
\subjclass{14D20, 14C34}

\address{School of Mathematics, Tata Institute of
Fundamental Research, Mumbai 400 005 (India)}

\email{indranil,tomas@math.tifr.res.in}

\begin{document}

\begin{abstract}
Let $X$ be a smooth projective curve over $\CC$, and let
$\SM^{n,\xi}_X$ be the moduli space of stable Higgs bundles on 
$X$ (with genus $g>1$), 
with rank $n$ and fixed determinant $\xi$, with $n$ and
$\deg(\xi)$ coprime. Let $X'$ and $\xi'$ be another such curve
and line bundle. We prove that
if $\SM^{n,\xi}_X$ and $\SM^{n,\xi'}_{X'}$ are isomorphic as
algebraic varieties, then $X$ and $X'$ are isomorphic.
\end{abstract}

\maketitle

\section{Introduction}

The classical Torelli theorem says that the Jacobian 
$J(X)$ of a curve, together with the polarization given
by the theta divisor, determines the curve $X$, i.e.,
if $J(X)$ and $J(X')$ are isomorphic as polarized varieties,
then $X$ is isomorphic to $X'$.

A similar result holds for $\SS\SU^{n,\xi}_X$, the moduli space of
stable vector bundles on $X$ with rank $n$ and fixed determinant
$\xi$, with $x$ and $\deg(\xi)$ coprime (assuming $g>1$). Namely,
the isomorphism class of $\SS\SU^{n,\xi}_X$ determines the
isomorphism class of $X$ (note that this moduli space has
a unique generator of polarization).
It was proved by 
Mumford and Newstead \cite{MN} for $n=2$, and later extended to any rank
by Narasimhan and Ramanan (\cite{NR} and \cite{NR2}). 
They consider the intermediate Jacobian 
associated to the third cohomology $H^3(\SS\SU^{n,\xi}_X)$.
It has a canonical polarization defined by the
positive generator of $\text{Pic}(\SS\SU^{n,\xi}_X)$.
They show that this canonically 
polarized intermediate Jacobian is isomorphic 
(as a polarized variety) to the 
Jacobian $J(X)$ of the curve with the polarization given
by the theta divisor, and then the result 
follows from the classical Torelli theorem.

In this paper we consider the same question for the moduli
space of Higgs bundles, with fixed determinant of degree coprime to
the rank, and traceless Higgs field.

Let $X$ be a connected smooth projective curve over $\CC$ (we will
assume that its genus $g>1$).
A Higgs bundle on $X$ is a pair $(E,\varphi)$, where $E$ is 
a vector bundle on $X$ and 
$$
\varphi: E\too E\otimes K_X
$$
is a morphism, called the Higgs field.
A subsheaf $F$ of $E$ is called $\varphi$-invariant if
$\varphi$ maps $F$ to $F\otimes K_X \subset E\otimes K_X$.
We say that a Higgs bundle is stable (respectively, semistable)
if for all $\varphi$-invariant proper subsheaves $F$ of $E$,
$$
\frac{\deg(F)}{\rk(F)} < \frac{\deg(E)}{\rk(E)} \quad 
(\text{respectively,} \; \leq).
$$
Denote by $\SM^{n,d}_X$ the moduli space of 
semistable Higgs sheaves with $\rk(E)=n$ and $\deg(E)=d$.
Note that the trace of $\varphi$ is a section of $K_X$.

We denote by $\SM^{n,\xi}_X$ the moduli space of semistable
Higgs bundles with $\tr(\varphi)=0$, $\rk(E)=r$, and
$\det(E)\isom \xi$, where $\xi$ is a fixed line bundle on 
$X$. Both $\SM^{n,d}_X$ and $\SM^{n,\xi}$ are irreducible and
if $r$ and $d$ (respectively, $\deg(\xi)$) are coprime, then
$\SM^{n,d}$ (respectively, $\SM^{n,\xi}$) is smooth.

The main result of this paper is the following theorem.

\begin{theorem}
\label{torelli}
Let $n>0$ and $d$ be two coprime integers.
Let $X$ and $X'$ be two smooth projective
curves with $g\geq 2$, 
and let $\xi$ and $\xi'$ be line bundles on $X$ and $X'$,
with $\deg(\xi)=\deg(\xi')=d$. 
If there is an isomorphism of algebraic varieties
$$
\SM^{n,\xi}_X \;\isom\; \SM^{n,\xi'}_{X'},
$$
then $X$ and $X'$ are isomorphic.
\end{theorem}

\noindent\textbf{Outline of the proof.}
Recall that there is a surjective morphism, called the
\textit{Hitchin map}, 
from the moduli space of Higgs bundles to a vector space
of dimension $(n^2-1)(g-1)$, called the \textit{Hitchin space}. 
The fiber over the origin is called
the nilpotent cone. It has several irreducible components,
and one of them is isomorphic to $\SS\SU^{n,\xi}$. In fact, it
is the only irreducible component
of the nilpotent cone that \textit{does not} admit a nontrivial
$\CC^*$ action
(cf. Proposition \ref{action}). Hence, if we were given 
$\SM^{n,\xi}_X$ together with the Hitchin map, we would
recover $\SS\SU^{n,\xi}_X$, and using the Torelli theorem
for vector bundles, we would also recover $X$.

The problem, of course, is that we are given only the 
isomorphism class of $\SM^{n,\xi}_X$,
but not the Hitchin map.

We start with an algebraic variety $Y$, isomorphic to
$\SM^{n,\xi}_X$. Look at the natural map
$$
m:Y \too \Spec \Gamma(Y),
$$
where $\Gamma(Y)$ is the ring of global functions. 
Since the
fibers of the Hitchin map are complete, it turns out that
this map is isomorphic to the Hitchin map.
More precisely, by Lemma \ref{functions}, 
$\Spec \Gamma(Y)\isom \AA^{(n^2-1)(g-1)}$, and there is
a commutative diagram like (\ref{isohitchin}).

This gives the Hitchin map, but only up to automorphism (as an
algebraic variety)
of $\AA^{(n^2-1)(g-1)}$. More precisely, we have recovered
the Hitchin fibration, but not the Hitchin map.
We have to find which point
of $\AA^{(n^2-1)(g-1)}$ corresponds to 
the origin of the Hitchin space.

Recall that the moduli space of Higgs bundles has a
$\CC^*$ action given by sending $(E,\varphi)$ to $(E,t\varphi)$.
This action descends to an action on the Hitchin space,
and the origin has the property that it is the only fixed
point of this descended action. In Proposition \ref{fixed}
we prove that this property characterizes the origin.
More precisely,
if $g$ is a $\CC^*$ action on the Hitchin space, having
exactly one fixed point, and admitting a lift to  
$\SM^{n,\xi}$, then this fixed point is the origin.
Hence, if $g$ is such an action on $\AA^{(n^2-1)(g-1)}$,
its fixed point is precisely the point corresponding
to the origin of the Hitchin space, the fiber of
$m$ over this point is isomorphic to the nilpotent
cone, and by the previous observations, the theorem is
proved.

Now we will give an outline of the proof of Proposition \ref{fixed}.
Let $s$ be a point in the Hitchin space such that the 
corresponding spectral curve $X_s$ is smooth. Then we have
a Kodaira-Spencer map $u_0$ from the tangent space at $s$ to the space of 
deformations of $X_s$. The kernel of this map is described in
Proposition \ref{kernel}. In the case
of $\SM^{n,\xi}_X$, this kernel has dimension one.
The direction defined by this kernel is precisely the direction
given by the standard $\CC^*$ action. 
The fiber over $s$ is isomorphic to the Prym variety $P_s$
(cf. (\ref{prym})), hence we also have a Kodaira-Spencer map $v_0$
from the tangent space at $s$ to the space of deformations
of $P_s$. It can be shown that a nontrivial deformation of the
spectral curve produces a nontrivial deformation of $P_s$,
and hence the kernel of the homomorphism $v_0$ coincides with
the kernel of $u_0$.

If $g$ is a $\CC^*$ action on $\AA^{(n^2-1)(g-1)}$ that lifts to
$\SM^{n,\xi}_X$, then
the tangent vector at $s$ defined by the action $g$ 
is in the kernel of the Kodaira-Spencer map $u_0$.
This implies that the orbit $g(\CC^*,s)$ is contained in
the orbit of the standard action, and then the origin of the
Hitchin space is in the closure of $g(\CC^*,s)$. Now,
a point in the closure of an orbit is a fixed point,
hence if $g$ has exactly one fixed point, this fixed
point has to be the origin. This proves Proposition \ref{fixed}. 

\section{Preliminaries}

In this section we will recall some general facts about
Higgs bundles (see, for instance, \cite{Hi} and \cite{BNR}). 
Define the Hitchin space
$$
\SH = H^0(K^{}_X) \oplus H^0(K_X^2) \oplus \cdots \oplus H^0(K_X^n).
$$
Its dimension is $n^2(g-1)+1$.

Let $S$ be the total space of the line bundle $K_X$, and
let $p$ be the projection
$$
p:S=\underline\Spec (\Sym^\bullet K_X^{-1})\too X.
$$
Given $(s)=(s_1,\ldots, s_n)\in \SH$, we define the spectral
curve $X_s$ in $S$ as the zero scheme of the following section
of $p^*K_X^n$
$$
f=x^n + \wt s_1 x^{n-1} + \wt s_2 x^{n-2} + \cdots + \wt s_n,
$$
where $\wt s_i = p^*s_i$, and $x\in H^0(S,p^*K_X)$ is the
tautological section. In other words, 
\begin{equation}
\label{xsspec}
X_s = \underline\Spec(\Sym^\bullet(K_X^{-1})/\SI)
\end{equation}
where $\SI$ is the ideal sheaf generated by the image of the 
homomorphism
$$
\begin{array}{rcl}
K_X^{-n} &\too &\Sym^\bullet(K_X^{-1})\\
\rule{0cm}{0.5cm}\alpha & \longmapsto & \alpha\sum_{i=0}^n s_i
\end{array}
$$
where we take $s_0=1$. Note that the projection $\pi:X_s\to X$,
which is the restriction of $p$,
has degree $n$. From (\ref{xsspec}) we obtain the following
isomorphism
\begin{equation}
\label{isopio}
\pi_*\SO^{}_{X_s} = \SO^{}_X \oplus K_X^{-1} \oplus K_X^{-2} 
\oplus \cdots \oplus K_X^{-(n-1)}
\end{equation}

There is a dense open set $U$ in $\SH$ (respectively, $U_0$ in $\SH_0$)
such that the spectral curve $X_s$ is smooth for $s\in U$
(respectively, $U_0$). Let $X_s$ be such a curve, and consider the
short exact sequence
$$
0 \too T_{X_s} \too T_S|_{X_s} \too N_{X_s/S} \too 0.
$$
Since $S$ is the cotangent space of $C$, it is a holomorphic 
symplectic variety,
and then $\det(T_S)\isom \SO_S$. Hence, we have
\begin{equation}
\label{iso0}
N_{X_s/S}\isom K_{X_s}.
\end{equation}
On the other hand,
\begin{equation}
\label{iso1}
N^{}_{X_s/S}\isom \SO(X_s)|^{}_{X_s} = p^*K_X^n|^{}_{X_s} = \pi^* K_X^n,
\end{equation}
and then the ramification line bundle of the projection 
$\pi:X_s\too X$ is
$$
\SO(R)=K^{}_{X_s} \otimes \pi^*K_X^{-1} = \pi^*K_X^{n-1}.
$$
The section
\begin{equation}\label{ramification}
\frac{\partial f}{\partial x}= n x^{n-1} + (n-1)\wt s_1 x^{n-2} +
\cdots + \wt s_{n-1} \;\in\; H^0(S,p^*K_X^{n-1}),
\end{equation}
when restricted to $X_s$, gives a section of $\SO(R)$, and its
scheme of zeroes is exactly the ramification divisor $R$.

Given a Higgs bundle $(E,\varphi)$, the formula
$$
\det(x\cdot \id-\varphi)=x^n+s_1 x^{n-1} + s_2 x^{n-2} +\cdots + s_n
$$
defines sections $s_i\in H^0(K_X^i)$, and this defines
the Hitchin map
$$
h: \SM^{n,d} \too \SH.
$$
The dimension of $\SM^{n,d}$ is $2n^2(g-1)+2$, and the 
fibers of $h$ are equidimensional projective schemes 
of dimension $n^2(g-1)+1$. If $X_s$ is smooth, then 
the fiber $h^{-1}(s)$ is isomorphic to the Jacobian $J(X_s)$
of the spectral curve \cite{Hi}, \cite{BNR}.

We can restrict this map to the subscheme $\SM^{n,\xi}\,\subset\,
\SM^{n,d}$ corresponding
to fixed determinant and traceless $\varphi$. Since
$s_1=\tr(\varphi)$, the image is actually in the traceless
Hitchin space 
\begin{equation}
\label{h0}
h_0: \SM^{n,\xi} \too \SH_0 = \bigoplus_{i=2}^n H^0(K_X^i).
\end{equation}
We have
\begin{eqnarray*}
\dim(\SM^{n,\xi})&=&2(n^2-1)(g-1),\\
\dim(\SH_0)&=&(n^2-1)(g-1),
\end{eqnarray*}
and the fibers of $h_0$ are equidimensional projective schemes 
of dimension $(n^2-1)(g-1)$. If the spectral curve $X_s$ is
smooth, then the fiber $h_0^{-1}(s)$ is the Prym variety
\begin{equation}
\label{prym}
P_s =\big\{ L\in \Pic(X_s): \det(\pi_*L)\isom \xi \big\},
\end{equation}
where $\pi$ is the obvious projection of $X_s$ to $X$.

\begin{lemma}
\label{covering}
Let $s\in \SH_0$ such that $X_s$ is smooth. Then
the morphism
\begin{eqnarray*}
\alpha: P_s \times J(X) & \too & J(X_s) \\
(L_1,L_2) & \longmapsto & L_1 \otimes \pi^* L_2
\end{eqnarray*}
is an unramified covering of degree $n^{2g}$.
\end{lemma}

\begin{proof}
First note that $\dim(J(X)\times P_s)=\dim (J(X_s))$.
If $L_1\otimes \pi^*L_2 \isom L'_1\otimes \pi^*L'_2$, 
then 
\begin{eqnarray*}
(\pi_* L_1) \otimes L_2 &\isom& (\pi_* L'_1) \otimes L'_2, \\
 \xi \otimes L_2{}^n &\isom&  \xi \otimes L'_2{}^n, 
\end{eqnarray*}
and hence $L'_2\isom L_2 \otimes \zeta$, for some $\zeta$
in  $J(X)[n]$  (the $n$-torsion subgroup of $J(X)$), 
and $L_1\isom L'_1
\otimes \pi^*(\zeta)$. Hence, the fiber of $\alpha$ is
isomorphic to $J(X)[n]$, and the lemma follows.
\end{proof}

\begin{lemma}
\label{functions}
The Hitchin map induces an isomorphism between the rings
of global functions
$$
\Gamma(\SM^{n,\xi})\;\isom\; \Gamma(\SH_0)
\;\isom\; \CC[y^{}_1,y^{}_2,\ldots, y^{}_{(n^2-1)(g-1)}].
$$
The last isomorphism follows from $\SH_0\isom \AA^{(n^2-1)(g-1)}$.
\end{lemma}

\begin{proof}
Since the Hitchin map is surjective, it gives an inclusion
$$
\Gamma(\SM^{n,\xi})\;\supset\;\Gamma(\SH_0).
$$
Since the fibers of $h_0$ are projective, any function on 
$\SM^{n,\xi}$ is constant on these fibers, hence it factors
through the Hitchin space, and hence we have an equality.
\end{proof}

\section{The nilpotent cone}

The fiber of the Hitchin map $h_0$ over the origin $s=0$ is called
the nilpotent cone. It is a reducible scheme of dimension
$(n^2-1)(g-1)$. Note that $(E,\varphi)$ is in the nilpotent
cone if and only if $\varphi$ is a nilpotent endomorphism.

\begin{proposition}
\label{action}
The nilpotent cone has a unique irreducible component that
does not admit a nontrivial $\CC^*$ action, and this is isomorphic
to $\SS\SU^{n,\xi}$, the moduli space of semistable vector bundles
on $X$ of rank $n$ and fixed determinant $\xi$.
\end{proposition}

Before proving this proposition, we need the following lemma
(recall that we are assuming that $r$ and $d$ are coprime, and hence
semistable implies stable).

\begin{lemma}[Lemma 11.9 in \cite{SiII}]
Let $(E,\varphi)$ be a Higgs bundle in the nilpotent cone,
with $\varphi\neq 0$. Consider the standard $\CC^*$ action
sending $(E,\varphi)$ to $(E,t\varphi)$. Assume that $(E,\varphi)$
is a fixed point, i.e., for every $t$ there is an isomorphism
with $(E,t\varphi)$.
Then there is another Higgs bundle $(F,\psi)$ in the nilpotent
cone, not isomorphic to $(E,\varphi)$, such that 
$\lim_{t\to \infty}(F,t\psi)=(E,\varphi)$.
\end{lemma}

\begin{proof}
This lemma is stated in \cite{SiII} with the extra assumption 
that $\deg(E)=0$ (because Simpson is interested in representations 
of the fundamental group), but the proof actually works for
Higgs bundles of any degree. For convenience of the reader,
we will give the necessary details.

Since $(E,\varphi)$ is a fixed point, by \cite[Lemma 4.1]{Si}
we know that it is of the form
$$
E = \bigoplus_{i=0}^m E_p
$$
with $\varphi$ sending $E_i$ to $E_{i-1}\otimes K_X$.
Since $\varphi\neq 0$, we have $m>0$ and $E_0$ and $E_m$ are
nontrivial. 
Since $r$ and $d$ are coprime, $(E,\varphi)$ is
stable. Then, since $E_0$ is $\varphi$-invariant, 
the stability condition implies 
$$
\frac{\deg(E_0)}{\rk(E_0)}<\frac{d}{r}
$$
The subbundle $\oplus_{i=1}^m E_p$ is also $\varphi$-invariant,
and the stability condition implies
$$
\frac{d}{r}<\frac{\deg(E_m)}{\rk(E_m)}.
$$
Combining both inequalities, we obtain $\deg\Hom(E_m,E_0)<0$,
and Riemann-Roch implies that $\Ext^1(E_m,E_0)$ is nontrivial.
Furthermore, if $A\subset E_m$ is the $\beta$-subbundle of 
$E_m$ (i.e., the first term in the Harder-Narasimhan filtration), 
then the
slope of $A$ is bigger or equal than the slope of $E_m$,
so we still have 
$$
 \Ext^1(E_m,E_0)\;\surj\; \Ext^1(A,E_0) \neq 0.
$$
Let $\eta$ be a nonzero element of $\Ext^1(E_m,E_0)$ that
maps to a nonzero element in
$\Ext^1(A,E_0)$. For each $t\in \CC^*$, let
$M_t$ be the extension given by $t^m \eta$
$$
0 \too E_0 \too M_t \too E_m \too 0.
$$
The Higgs bundle $(F_t,\psi_t)$ is defined taking 
$$
F_t = M_t \oplus \bigoplus_{0<p<m} E_p
$$
and the Higgs field $\psi$ is given by
\begin{eqnarray*}
& E_i \stackrel{\varphi}\too E_{i-1}\otimes K_X , & 1<i<m\\
&E_1 \stackrel{\varphi}\too E_0\otimes K_X \too M_t \otimes K_X& \\
&M_t \too E_m \stackrel{\varphi}\too E_{m-1}\otimes K_X&
\end{eqnarray*}
By construction, we see that $\psi_t$ is nilpotent.
We define $(F,\psi)=(F_t,\psi_t)$.
The rest of the proof is identical to the proof given in 
\cite{SiII}, so we only give a sketch.
First one checks that $(F,t^{-1}\psi)$ is isomorphic to 
$(F_t,\psi_t)$, and that these Higgs bundles are stable,
hence they are in the nilpotent cone. Furthermore,  
$$
\lim_{t\to \infty}(F,t\psi)=(E,\varphi).
$$
Finally Simpson checks that $(F,\psi)$ is not isomorphic to 
$(E,\varphi)$ (here is where we use the fact that 
$\eta$ has nonzero image in $\Ext^1(A,E_0)$).
\end{proof}

Now we can prove the proposition.

\begin{proof}[Proof of Proposition \ref{action}]
The map that sends a vector bundle $E$ to the Higgs bundle $(E,0)$
defines an inclusion of $\SS\SU^{n,\xi}$ in the nilpotent cone.
Since they have the same dimension, this gives one component
of the nilpotent cone. It does not have a nontrivial $\CC^*$ action,
because if it had, then it would produce a nontrivial section
of the tangent bundle of $\SS\SU^{n,\xi}$, but this is known
to have no nontrivial sections (\cite{NR} and \cite{NR2}).

In the rest of the components we have the $\CC^*$ action given
by $(E,\varphi)\mapsto (E,t\varphi)$. To show that this 
action is nontrivial, we use the previous lemma:
Let $(E,\varphi)$ be
a fixed point (with $\varphi\neq 0$), and let $(F,\psi)$
be the Higgs bundle given by Lemma 3.2. 
Since the limit of the
action moves $(F,\psi)$ to $(E,\varphi)$, they are
in the same irreducible component of the nilpotent cone, and 
since both bundles
are stable and nonisomorphic, they correspond to different
points in the moduli space, and hence the $\CC^*$ action
is nontrivial in this component.
\end{proof}

\section{Kodaira-Spencer map}

The projection $\pi:X_s \too X$ has degree $n$. 
The spectral curve construction gives a bijection between 
points in $\SH$ and projective curves on $S$ such that the
restriction of the projection has degree $n$.

Let $s\in \SH$ be a point such that the corresponding spectral
curve $X_s$ is smooth. There is a morphism from the tangent
space $T_s\SH\isom \SH$ to the space $H^1(X_s, T_{X_s})$ of
all infinitesimal deformations
of $X_s$, the Kodaira-Spencer map $u$.

We can also define the Kodaira-Spencer map $u_0$ for the 
traceless Hitchin space $\SH_0$ (cf. (\ref{h0})). It is 
obtained by restricting $u$ to $T_s\SH_0 \subset T_s\SH$.

\begin{lemma}
\label{th}
If $X_s$ is smooth, then there are natural isomorphisms
$$
H^0(X_s,N_{X_s/S})\;\isom\; H^0(X_s,\pi^*K_X^i) \;\isom\; T_s\SH .
$$
\end{lemma}

\begin{proof}
Using the isomorphisms (\ref{iso1}) and (\ref{isopio}) and the projection 
formula, we have
\begin{eqnarray*}
&H^0(X_s,N_{X_s/S}) = H^0(X_s,\pi^*K_X^n) =
H^0(X,K_X^n\otimes \pi_*\SO_{X_s}) =&\\
&= H^0(X,\oplus_{i=1}^n K_X^n) = \SH \isom T_s\SH&
\end{eqnarray*}

\end{proof}

The objective of this section is to calculate the kernels
of $u$ and $u_0$. 
There are some elements in $H^0(X_s,N_{X_s/S})$ that are
clearly in the kernel. For instance, let 
$\lambda\in \CC\isom H^0(X, \SO_X)$,
and denote a point in $X_s\subset S$ by $(\omega,x)$, where 
$\omega$ is a point in $X$ and $x$  is a coordinate
in the fiber of $S$ over $\omega$. Then the deformation sending 
$(\omega,x)$ to $(\omega,e^\lambda x)$ clearly does not change the
isomorphism class of $X_s$. In fact, this is the deformation 
produced by the standard $\CC^*$ action, and it is clearly
in the kernel of the Kodaira-Spencer map $u_0$.

Furthermore, for any $\alpha\in H^0(X,K_X)$, sending $(\omega,
x)$ to $(\omega,x+\alpha(\omega))$
also preserves the isomorphism class of $X_s$. The deformations
defined in this way do not preserve the condition
$0=\tr(\varphi)$ $(=s_1)$, and hence they are in the kernel of $u$,
but not in the domain of $u_0$. The following proposition
says that these two constructions describe the kernels.

\begin{proposition}
\label{kernel}
There is an exact sequence
$$
0 \too H^0(X,K_X \oplus \SO_X) \too T_s\SH 
\stackrel{u}{\too} H^1(X_s,T_{X_s}),
$$
hence the kernel of the Kodaira-Spencer map $u$ has dimension $g+1$.
If we fix the determinant, we have an exact sequence
$$
0 \too H^0(X,\SO_X) \too  T_s\SH_0 
\stackrel{u_0}{\too} H^1(X_s,T_{X_s}),
$$
and hence the kernel of the restricted Kodaira-Spencer map
$u_0$ has dimension $1$.
\end{proposition}

\begin{proof}
Consider the following diagram, constructed using (\ref{iso0})
and (\ref{iso1})
$$
\xymatrix{
  & & {p^*T_{X}|_{X_s}}  & & \\
{0} \ar[r] &  {T_{X_s}} \ar[r] & 
{T_S|_{X_s}} \ar[rr]\ar@{>>}[u] && {N_{X_s/S}}
\ar[r] & 0 \\
 & & {\rule{20pt}{0pt}T_p|_{X_s}\isom \pi^*K_X} \ar@{^{(}->}[u] 
\ar[rr]^{\otimes\frac{\partial f}{\partial x}\big|^{}_{X_s}} 
&& {\pi^*K_X^n\isom K^{}_{X_s}\rule{20pt}{0pt}} \ar@{=}[u] & 
}
$$
where $T_p$ denotes the relative tangent bundle for the
projection $p$.
Note that the diagram is well defined, since $\frac{\partial f}{\partial
x}\big|_{X_s}$ is a section of $\pi^*K_X^{n-1}\isom \SO(R)$ (cf. 
(\ref{ramification})).
The diagram is commutative because the zero scheme
of the two morphisms between the line bundles $T_p|_{X_s}$ and 
$N_{X_s/S}$
are the same, namely the ramification divisor, hence the maps
differ by a scalar, but this scalar is absorbed in the 
isomorphism between $K_{X_s}$ and $N_{X_s/S}$.

Since the tangent line bundle $T_X$ has negative degree,
$H^0(X_s, p^*T_{X}|_{X_s})=0$. Therefore,
the previous diagram gives
$$
\xymatrix{
{0} \ar[r] & {H^0(T_S|_{X_s})} \ar[r]^{I} &
{H^0(N_{X_s/S})} \ar[r]  & {H^1(T_{X_s})} \ar@{=}[dd] \\
 & H^0(X_s,\pi^*K_X) \ar[u]^{\isom} \ar[r]^{I'}   & 
H^0(X_s,\pi^*K_X^n) \ar[u]^{\isom}_{\beta} & \\
{0} \ar[r] & H^0(X,K_X\oplus \SO_X)\ar[u]^{\isom}_{\gamma} \ar[r]^-{I''} & 
{T_s\SH} \ar[u]^{\isom}_{\alpha} \ar[r]^{u}& {H^1(T_{X_s})} 
}
$$
where $\alpha$ and $\beta$ are the isomorphism of Lemma \ref{th},
the isomorphism $\gamma$ is given by
\begin{eqnarray*}
&&H^0(X_s,\pi^*K_X) = H^0(X,K_X\otimes \pi_*\SO_{X_s})=\\
&&\rule{1cm}{0cm} =  H^0(X,\oplus_{i=1-n}^1 K_X^i) =
H^0(X, K_X \oplus \SO_X),
\end{eqnarray*}
$I'=H^0(\otimes \frac{\partial f}{\partial x}\big|_{X_s})$, 
$I''$ is defined by composition, 
and $u$ is the Kodaira-Spencer map. The bottom row is the
first exact sequence in the statement of the proposition.

Now we will restrict the Kodaira-Spencer map $u$ to $T_s\SH_0$.
To do this, we need an explicit description of the morphism $I'$.

Using the isomorphism $H^0(X_s,\pi^*K_X^n)\isom 
H^0(X,\oplus_{i=1}^n K_X^n)$ (cf. proof of Lemma \ref{th}),
an element of this group is written as
$$
\wt a_0 x^n + \wt a_1 x^{n-1} + \cdots + \wt a_n
$$
with $a_i\in H^0(X, K_X^i)$ and $\wt a_i=\pi^* a_i$.

On the other hand, using the isomorphism $\gamma$, 
an element of this group will be written as
$$
\wt b_1 + \wt b_0 x
$$
with $b_i\in H^0(X, K_X^i)$ and $\wt b_i=\pi^* b_i$.

Since $I'$ comes from multiplication with 
$\frac{\partial f}{\partial x}\big|_{X_s}$, a short calculation
using $f|_{X_s}=0$ gives
\begin{equation}
\label{exfor}
I'(\wt b_1 + \wt b_0 x) \;=\;
\sum_{i=1}^{n} 
\big( (n-i+1) \wt s_{i-1}\wt b_1 - i \wt s_i \wt b_0\big)x^{n-i}
\end{equation}

The subspace $\SH_0\subset \SH$ is the zero locus of the
trace map sending $(s_1,\ldots,s_n)$ to $s_1$. Then we
have a commutative diagram
$$
\xymatrix{
0 \ar[r] & {T_s\SH_0} \ar[r] & {T_s\SH} \ar[r]^{d(\tr)} & 
H^0(X,K_X) \ar[r] & 0 \\
0 \ar[r] & {H^0(\oplus_{i=2}^n K_X^i)} \ar[r]\ar[u]^{\isom} & 
{H^0(\oplus_{i=1}^n K_X^i)} \ar[r]^{p_1} \ar[u]^{\isom} & 
H^0(X,K_X) \ar[r] \ar@{=}[u] & 0 
}
$$
where $p_1$ is projection to the first summand.

Now, if $s\in \SH_0$, then $s_1=0$, and using the explicit
formula (\ref{exfor}) for $I'$, we obtain 
$$
(d(\tr) \circ I'') (\wt b_1 + \wt b_0 x) \;=\; n\wt b_1,
$$
and hence the following diagram is commutative
$$
\xymatrix{
& 0 \ar[d] & 0 \ar[d]  \\
0 \ar[r] & {H^0(\SO_X)} \ar[r] \ar[d] & {T_s\SH_0} \ar[r]^-{u_0}\ar[d]
& H^1(X_s,T_{X_s})  \ar@{=}[d]  \\
0 \ar[r] & {H^0(K_X \oplus \SO_X)} \ar[r]^-{I''} \ar[d]^{q} & 
{T_s\SH} \ar[r]^-{u}\ar[d]^{d(\tr)} & {H^1(X_s,T_{X_s})}  \\
 & {H^0(K_X)} \ar[d]\ar@{=}[r] & {H^0(K_X)}\ar[d] &  \\
 & 0 & 0 & \\
}
$$
where $q$ is projection to the first summand followed by
multiplication by $n$. The top row is the
second exact sequence in the statement of the proposition.
\end{proof}

\section{Torelli theorem}

Recall that there is a $\CC^*$ action on $\SM^{n,\xi}$,
sending $(E,\varphi)$ to $(E,t\varphi)$. This action 
descends to give an action on $\SH_0$, whose only fixed
point is the origin $s_i=0$.

\begin{proposition}
\label{fixed}
Let $g:\CC^*\times \SH_0 \too \SH_0$ be an action, 
having exactly one fixed,
and admitting a lift to $\SM^{n,\xi}$.
Then this fixed point is the origin $s_i=0$.
\end{proposition}

\begin{proof}
As we have seen, the standard action has this property.
Now let $s\in \SH_0$ be a point giving a smooth spectral 
curve $X_s$. The tangent vector at $s$ defined by the standard action
is contained in the kernel of the Kodaira-Spencer map $u_0$, because
the standard action does not change the isomorphism class of the 
spectral curve. Now we will show that the tangent vector
defined by any action that lifts to $\SM^{n,\xi}$ is also
in the kernel of the Kodaira-Spencer map.

Denote $J=J(X)$ and $J_s=J(X_s)$.
In general, if $M_1\to M_2$ is a covering map between differentiable
manifolds, a holomorphic structure on $M_2$ induces a holomorphic
structure on $M_1$. Hence, using the morphism $\alpha$ of
Lemma \ref{covering}, a deformation of $J_s$ gives a deformation of
$J\times P_s$. The map $\epsilon$ between the deformation spaces is 
the composition 
$$
\xymatrix{
{\epsilon: H^1(J_s, T_{J_s}) } \ar@{^{(}->}[r] &
{H^1(J_s, T_{J_s}\otimes \alpha_*\SO_{J\times P_s})} \ar[r]^-{\isom}&
 {H^1(J\times P_s, T_{J\times P_s})}
}
$$
where the inclusion is induced by the adjunction map 
$\SO_{J_s}\to \alpha_* \SO_{J\times P_s}$ and the isomorphism is given
by the projection formula and the isomorphism $\alpha^* T_{J_s}\isom
T_{J\times P_s}$.

We have several varieties associated to a point $s\in \SH_0$, namely
the spectral curve $\pi:X_s\too X$, the Jacobian $J_s$ and the 
Prym variety $P_s$. Hence the $g$ action produces deformations of all
these objects
\begin{eqnarray*}
\eta_1 & \in & H^1(X_s,T_{X_s}) \\
\eta_2 & \in & H^1(J_s,T_{J_s}) \\
\eta_3 & \in & H^1(P_s,T_{P_s})
\end{eqnarray*}
We will study the relationship between $\eta_2$ and $\eta_3$.
Using the projection formula and the Leray spectral sequence 
for the projections $q_1$ and 
$q_2$ from $J\times P_s$ to $J$ and $P_s$, we obtain that
the space of deformations of $J\times P_s$ is
\begin{eqnarray}
\label{summands}
H^1(J\times P_s,T_{J\times P_s}) = H^1(J\times P_s,q_1^*T_J) 
\oplus H^1(J\times P_s,q_2^*T_{P_s}) = \\
= H^1(T_J) \oplus \Big( H^0(T_J)\otimes H^1(\SO_P) \Big)
\oplus \Big( H^0(T_P)\otimes H^1(\SO_J) \Big) \oplus
H^1(T_{P_s}) 
\nonumber
\end{eqnarray}
The first and last terms correspond to deformations of $J$ and $P_s$.
To understand the second term, view $J\times P_s \to P_s$ as a 
trivial fiber bundle. The deformations of this as a fiber bundle
(i.e., keeping the fiber and the base fixed) are parametrized by
$$
H^1(P,\SO_P\otimes\Lie(\Aut(J))=H^1(P,\SO_P\otimes H^0(J,T_J))=
H^1(P,\SO_P) \otimes H^0(J,T_J),
$$
i.e., the second term. Analogously, the third term corresponds
to deformations of $J\times P_s \to J$ as a fiber bundle.

The class $\eta_3$ lies in the fourth summand 
$H^1(P_s,T_{P_s})$. Indeed, if we move $s$ to a nearby point $s'$,
the fiber $J\times P_s$ is deformed to $J\times P_{s'}$, hence
the components of $\eta_3$ in the first three summands of 
(\ref{summands}) have to be zero.
Hence $\epsilon(\eta_2)$ lies in $H^1(P_s,T_{P_s})$, and is 
equal to $\eta_3$.
$$
\epsilon(\eta_2)=\eta_3 \;\in\;  H^1(P_s,T_{P_s}) 
\;\subset\; H^1(J\times P_s,T_{J\times P_s}) \\
$$

A deformation of a curve produces a deformation of its 
Jacobian, hence there is a natural map
$$
H^1(X_s, T_{X_s}) \too H^1(J_s, T_{J_s}),
$$
and the infinitesimal version of the classical Torelli theorem says
that this map is injective. Clearly, in our situation $\eta_2$ is
the image of $\eta_1$ under this map. Combining the injectivity
of this map with the injectivity of $\epsilon$, we obtain that
if $\eta_3=0$, then $\eta_1=0$.

Since the action $g$ lifts, the fiber $P_s$ of the Hitchin map $h_0$ 
at $s$ is isomorphic to the fiber $P_{g(t,s)}$ at $g(t,s)$ 
for all $t\in \CC^*$. This implies that $\eta_3=0$, hence $\eta_1=0$.

This means that the tangent vector at $s$ defined by the action
$g$ is in the kernel of the Kodaira-Spencer map $u_0$, and since
by Proposition \ref{kernel} this has dimension $1$, we obtain that
the orbit $g(\CC^*,s)$ of $g$ through $s$ is included in the 
orbit of the standard
action through $s$. In particular, the origin is a limiting point
of the orbit $g(\CC^*,s)$, i.e., it is in the closure of the orbit,
but not in the orbit (note that the origin is not in the orbit,
because the fiber over the origin is not isomorphic to the fiber over
$s$). The limiting points of an orbit are fixed points of the 
action. Then the origin is a fixed point of the action, and 
by hypothesis is the only fixed point.
\end{proof}

Finally we prove the Torelli theorem for Higgs
bundles.

\begin{proof}[Proof of Theorem \ref{torelli}]
Let $Y$ be an algebraic variety isomorphic to
the moduli space $\SM^{n,\xi}_X$.
By Lemma \ref{functions}, choosing a set of generators of the
ring of global functions on $Y$ we obtain an isomorphism
$\Gamma(Y)\isom \CC[y^{}_1,y^{}_2,\ldots, y^{}_{(n^2-1)(g-1)}]$, 
and then the natural morphism
$Y \too \Spec \Gamma(Y)$ gives a morphism
$$
m:Y \too \AA^{(n^2-1)(g-1)}.
$$
Note that $m$ depends on the set of generators chosen: 
a different choice will give a morphism that differs by an automorphism
(as an algebraic variety) of $\AA^{(n^2-1)(g-1)}$.
Up to isomorphism, this is the Hitchin map. More precisely,
if $\alpha:Y \too \SM^{n,\xi}_X$ is an isomorphism, then there
is an isomorphism (as algebraic varieties) $\beta$ such that
the following diagram commutes
\begin{equation}
\label{isohitchin}
\xymatrix{
Y \ar[r]^{\alpha}_{\isom} \ar[d]_{m} & {\SM^{n,\xi}_X} \ar[d]^{h_0} \\
{\AA^{(n^2-1)(g-1)}} \ar[r]^-{\beta}_-{\isom} & \SH_0
}
\end{equation}
Let $g:\CC^*\times \AA^{(n^2-1)(g-1)} \too \AA^{(n^2-1)(g-1)}$
be a $\CC^*$ action with exactly one fixed point $y$, and such that
it admits a lift to $Y$. 

We know that such an action exists (using the standard $\CC^*$
action
on $\SH_0$ and the isomorphism $\beta$), and by Proposition
\ref{fixed}, $\beta(y)$ is the origin $s_i=0$ of $\SH_0$.
Indeed, if $\beta(y)$ were a different point, using 
$g$ and $\beta$ we would have an action contradicting
Proposition \ref{fixed}.

Then the fiber $m^{-1}(y)$ over $y$ is isomorphic to the nilpotent
cone. This has several irreducible components, but by
Proposition \ref{action}, only one of them (call it $Z$) 
does not admit nontrivial $\CC^*$ actions, and furthermore
$Z$ is isomorphic to $\SS\SU_X^{n,\xi}$, the moduli space of 
semistable vector bundles on $X$ of rank $n$ and fixed 
determinant $\xi$.

Now, if $\SM^{n,\xi}_X$ is isomorphic to $\SM^{n,\xi'}_{X'}$,
then $Z$ is isomorphic to $Z'$, hence
$\SS\SU_X^{n,\xi}$ is isomorphic to $\SS\SU_{X'}^{n,\xi'}$,
and by the Torelli theorem for the moduli space of vector
bundles, $X$ is isomorphic to $X'$.
\end{proof}


\begin{thebibliography}{EMGH}

\bibitem[BNR]{BNR}{A. Beauville, M.S. Narasimhan and S. Ramanan, }
\textit{Spectral curves and the generalised theta divisor,} 
J. Reine Angew. Math. \textbf{398} (1989), 169--179.

\bibitem[Hi]{Hi}{N. Hitchin, } \textit{Stable bundles and integrable
systems}, Duke Math. Jour. \textbf{54} (1987), 91--114.

\bibitem[MN]{MN}{D. Mumford and P.E. Newstead, }
\textit{Periods of a moduli space of bundles on curves, }
Amer. J. Math. \textbf{90} (1968), 1200--1208.

\bibitem[NR]{NR}{M.S. Narasimhan and S. Ramanan, }
\textit{Deformations of the moduli space of vector bundles
over an algebraic curve, } 
Ann. Math. (2) \textbf{101} (1975), 391--417. 

\bibitem[NR2]{NR2}{M.S. Narasimhan and S. Ramanan, }
\textit{Generalized Prym varieties as fixed points, }
Jour. Indian Math. Soc. \textbf{39} (1975), 1--19.

\bibitem[Si]{Si}{C. Simpson, }
\textit{Higgs bundles and local systems, }
Publ. Math. I.H.E.S. \textbf{75} (1992), 5--95.

\bibitem[SiII]{SiII}{C. Simpson, }
\textit{Moduli of representations of the fundamental group of a smooth
projective variety II, }
Publ. Math. I.H.E.S. \textbf{80} (1995), 5--79.

\end{thebibliography}
\end{document}